\documentclass[12pt,a4paper,reqno]{amsart}
\usepackage{amssymb,amsfonts,amsmath,amscd,amsthm,latexsym}
\usepackage{cite}
\usepackage{times}

\usepackage[bookmarks=false]{hyperref}
\hypersetup{%
    colorlinks=true,        % false: boxed links; true: colored links
    linkcolor=blue,          % color of internal links (change box color with linkbordercolor)
    citecolor=blue,         % color of links to bibliography
    urlcolor=blue           % color of external links
    }

\usepackage{url}
\newtheorem{theor}{Theorem}%[section]

\theoremstyle{definition}

\newtheorem{prop}[theor]{Proposition}%[section]
\newtheorem{deff}{Definition}%[section]

\newtheorem*{comment}{Comment}
\newtheorem*{observation}{Observation}

\newtheorem{open}{Open problem}

\newtheorem{ex}{Example}%[section]
\theoremstyle{remark}
\theoremstyle{definition}
\theoremstyle{definition}
\newcommand{\BBR}{\mathbb{R}}

\newcommand{\bx}{{\boldsymbol{x}}}

\newcommand{\veps}{\varepsilon}

\DeclareMathOperator{\Van}{Van}
\DeclareMathOperator{\Aut}{Aut}

\newcommand{\by}[1]{\textrm{{#1}}}
\newcommand{\jour}[1]{\textit{{#1}}}
\newcommand{\vol}[1]{\textbf{{#1}}}
\newcommand{\book}[1]{\textit{{#1}}}

\allowdisplaybreaks

\begin{document}

\title[Open problems about Jacobian determinants]{Kontsevich graphs act on Nambu\/--\/Poisson brackets,\\[3pt] VI. %Vanishing graphs and 
Open problems}

% TITLE OPTIONS:
% A- Vanishing graphs and open problems
% B- Why/How vanishing graphs vanish and other mysteries 
% C- An investigation into vanishing graphs
% D- Neatly vanishing graphs 
% E- How and why (do) graphs vanish
% F- 

%\runningheads{M.S.\ Jagoe Brown, A.V.\ Kiselev}{Open problems about Jacobian determinants}
% The short title can be whatever comes after "VI" in the main title %%%

%\begin{start}{%
\author[M.\,S.\,Jagoe Brown]{Mollie S.\ Jagoe Brown${}^{*}$}%{},
\thanks{${}^{*}$\:\textit{Present address}: Korteweg\/--\/de Vries Institute of Mathematics, University of Amsterdam, P.O.\ Box 94248, 1090\:GE Amsterdam, The Netherlands}
% The second argument connects author(s) with addresses

\author[A.\,V.\,Kiselev]{Arthemy V.\ Kiselev${}^{\S}$}%{1}
\thanks{${}^{\S}$\:\textit{Address}: Bernoulli Institute for Mathematics, Computer Science and 
Artificial Intelligence, University of Groningen, P.O.\,Box 407, 9700\:AK Groningen, The Netherlands}
% The second argument connects author(s) with addresses

%\email{m.s.jagoe.brown@gmail.com, a.v.kiselev@rug.nl}

% MSC codes: 
\subjclass[2010]{05C31, 05C90, 18G85, 53D55, 68R10} %Explanations:
% 05C31 - Graph polynomials; 05C90 - Applications of graph theory; 18G85 - Graph complexes and graph homology; 53D55 - Deformation quantization, star products (?); 68R10 - Graph theory (including graph drawing) in computer science (?)

\date{31 October 2025%
%Day Month Year (Insert date of submission)
}
% Insert date of submission
%}

\keywords{Casimir function,
differential\/-\/polynomial identity,
graph cocycle,
%graph calculus
Jacobian determinant,
Kontsevich's graph complex,
Levi\/-\/Civita symbol,
Nambu\/--\/Poisson bracket,
subgraph,
tensor\/-\/valued invariant of $GL(d)$ %-\/group 
action}

\begin{abstract}
Kontsevich's graphs from deformation quantisation allow encoding multi\/-\/vectors whose coefficients are differential\/-\/polynomial in components of Poisson brackets on finite\/-\/dimensional affine manifolds. The calculus of Kontsevich graphs can be made dimension\/-\/specific for the class of Nambu\/--\/Poisson brackets given by Jacobian determinants.
Using the Kontsevich\/--\/Nambu micro\/-\/graphs in dimensions $d\geqslant 2$,
we explore the open problem of (non)\/triviality for Kontsevich's tetrahedral graph cocycle action on the space of Nambu\/--\/Poisson brackets. We detect a conjecturally infinite new set of differential\/-\/polynomial identities for Jacobian determinants of arbitrary sizes~$d\times d$.%, and we indicate how they hold.
\end{abstract}
\maketitle

%\date{}

%\affil{Bernoulli Institute for Mathematics, Computer Science and Artificial Intelligence, University of Groningen, P.O. Box 407, 9700 AK Groningen, The Netherlands}

% Extremely rough first draft of new abstract:
% We found that vanishing graphs in dimension 3D are a crucial set to consider when expanding graphs to 4D in order to find a solution to the coboundary equation $Q^\gamma_d(P)=\llbracket P,\vec{X}^\gamma_d(P)\rrbracket$. Namely, just because a graph is vanishing in 3D does not mean that its descendants in 4D are too. Excluding vanishing graphs in a lower dimension can have a disastrous impact in higher dimensions. Here, we examine these extremely important elements in the study of Kontsevich graph deformations. We specifically propose a partial answer for why they vanish, and we detect how all of these vanishing 3D sunflower graphs vanish. 

\section{Introduction}\label{SecIntroduction}\noindent%
The Nambu\/-\/determinant Poisson bracket $\{f,g\}_d$ on~$\BBR^d\ni\bx=(x^1$,$\ldots$,$x^d)$, 
%is 
\begin{equation}\label{EqNambuBr}
\{f,g\}_d (\bx) = 
\varrho(\bx) \cdot \det \smash{\Bigl( }
\partial \bigl( f,g,a^1,\ldots,a^{d-2}\bigr) 
\big/ \partial \bigl( x^1, \ldots, x^d \bigr) \smash{\Bigr)}
,
\end{equation}
was introduced in~\cite{Nambu1973}; % in the context of quarks,
its Casimirs $a^i\in C^1(\BBR^d)$ satisfy $\{f,a^i\}_d\equiv 0$.
%%%%
%the contraction $\langle \varrho(\bx)\,\partial_{\bx}$, $f\otimes g\otimes a^1 \otimes \ldots \otimes a^{d-2} \rangle$ of a $d$-vector with a $d$-tuple of scalar functions.\footnote{\label{FootNambuFormula}
%In global coordinates $\bx=(x^1$,$\ldots$,$x^d)$ on~$\BBR^d$, the Nambu\/--\/Poisson bracket of $f,g\in C^1(\BBR^d)$ is
%\begin{equation}\label{EqNambuBr}
%\{f,g\}_d (\bx) = 
%\varrho(\bx) \cdot \det \Bigl( \partial \bigl( f,g,a^1,\ldots,a^{d-2}\bigr) 
%\big/ \partial \bigl( x^1, \ldots, x^d \bigr) \Bigr);
%\end{equation}
%the Casimirs $a^i\in C^1(\BBR^d)$ satisfy $\{f,a^i\}_d\equiv 0$.}
%%%
Following~\cite{Ascona96}, we study infinitesimal symmetries of the Jacobi identity, i.e.\ deformations of brackets which preserve their property to be Poisson (see~\cite{skew21,skew23,MSJB} and references therein).
By using `good' cocycles in the graph complex, Kontsevich revealed in~\cite{Ascona96} a class of such symmetries; the tetrahedron $\gamma_3$ is the smallest `good' %graph
cocycle.
The conjecture in~\cite{skew21} is that the $\gamma_3$-flow \emph{restricts} to the Nambu subset of Poisson brackets on~$\BBR^{d\geqslant 3}$.
Recall the standard problem in deformation theory: is the $\gamma_3$-flow \textit{non}\/tri\-vi\-al, i.e.\ does the change of Poisson bivector not amount to a change of coordinates along a vector field $\smash{\vec{X}^{\gamma_3}_d}$ on~$\BBR^d$\,?
Solution $\smash{\vec{X}^{\gamma_3}_{d=2}}$ was hinted in~\cite{Ascona96}, the field
$\smash{\vec{X}^{\gamma_3}_{d=3}}$ is known from~\cite{skew21}, and $\smash{\vec{X}^{\gamma_3}_{d=4}}$ was found in~\cite{MSJB}.
%%%
In this note %from the cycle I.--V. started in [I], cf.~[IV],
we sum up the surprising properties of graphs that encode solutions $\smash{\vec{X}^{\gamma_3}_{d\geqslant 2}(\varrho,\boldsymbol{a})}$ for the $\gamma_3$-flow of Nambu\/--\/Poisson brackets.\footnote{\label{FootOpenPrbNontriv}
It remains an open problem to find an example of \emph{nontrivial} graph cocycle action on a Poisson bracket, thus spreading it to a family of inequivalent structures not related by a coordinate change.
Beyond the tetrahedron $\gamma_3$, there are countably many `good' nontrivial graph cocycles: $\gamma_5$,\ $\gamma_7$,\ $\ldots$ (see~\cite{skew21,skew23} and references therein); they stem from the %generators of 
Grothendieck\/--\/Teichm\"uller Lie algebra~$\mathfrak{grt}$.%
}

We know the solutions $\smash{\vec{X}^{\gamma_3}_{d\leqslant 4}}$ separately in each dimension~$d$. There can be no universal formula of $\smash{\vec{X}^{\gamma_3}}$ --\,for all Poisson brackets in all dimensions\,-- for the nontrivial graph cocycle $\gamma_3$ of~\cite{Ascona96}; we do not know a universal formula of $\smash{\vec{X}^{\gamma_3}_d}$ --\,working in each %every 
$d\geqslant 3$\,-- for Nambu brackets~\eqref{EqNambuBr} in particular.\footnote{\label{FootReduceSizePrb}
To obtain the solution $\smash{\vec{X}^{\gamma_3}_{d=4}}$ over~$\BBR^4$, in~\cite{MSJB} we reduced the size of the problem circa %some 
300 times w.r.t.\ the worst\/-\/case scenario in~\cite{skew23}. This reduction effort continued in~\cite{IV} for steps $d\mapsto d+1$ in any %arbitrary 
dimension, given a solution $\smash{\vec{X}^{\gamma_3}_{d}}$ and some %certain 
`invisible' structures over~$\BBR^d$, which we now tackle.}
%%%
Studying the step $d\mapsto d+1$, we saw in~\cite{IV} that solution $\smash{\vec{X}^{\gamma_3}_{d=4}}$ is found \emph{economically} by knowing 
(\textit{i}) the combinatorics of formulas in $\smash{\vec{X}^{\gamma_3}_{d=2}}$ and $\smash{\vec{X}^{\gamma_3}_{d=3}}$
 %in terms of Nambu brackets
and (\textit{ii}) the structures that %which 
encode identically vanishing\footnote{\label{FootAnonceHowVanish}
Examples of this, not exhausting the full range, are given by formulas, differential\/-\/polynomial in the components of Poisson tensor, which equal minus themselves under a relabelling of summation indices. We stress the existence of other mechanisms for objects' vanishing.%; to reveal one of them is the problem which we presently address.
}
%%%
$1$-vectors built of whole Nambu brackets over~$\BBR^3$.
%%%
%In fact, 
These vanishing objects are tensor\/-\/valued invariants of $GL(d)$-%{} or affine group $\Aff(d)$ 
action; we study their collective behaviour under $d\mapsto d+1$ and internal construction of individual invariants. For them, we discover a likely infinite set of new identities, differential\/-\/polynomial w.r.t.\ the Casimirs $a^i$ 
in the Jacobian determinants and coefficients $\varrho(\bx)$ of $d$-vectors.
%%%
In this note we describe the procedure of (micro-)\/graph embeddings to obtain such identities and %we 
illustrate how they hold.

\noindent\textbf{Preliminaries.}\quad
Kontsevich's directed graphs from~\cite{Ascona96} encode polydifferential operators --\,in practice, multivectors\,-- made of copies of Poisson bi\/-\/vectors~$P$ as %the
building blocks: the subgraph of~$P$ is the wedge ${\leftarrow}\!{\bullet}\!{\rightarrow}$.
By definition, each digraph~$\Gamma$ (possibly, in an $\BBR$-linear combination) consists of the ordered set of $m$ sinks (where operator's arguments are placed) and $n$ internal vertices, the wedge tops; for each top, the pair of outgoing edges is ordered Left${}\prec{}$Right. Label the ordered sinks by $\mathsf{0}$,\ $\ldots$,\ $\mathsf{m-1}$ and label the wedge tops by $\mathsf{m}$,\ $\ldots$,\ $\mathsf{m+n-1}$. The encoding of~$\Gamma$ is the ordered list of $n$ ordered pairs\footnote{\label{FootSwapEdgesSign}
Swapping the order of two outgoing edges in a wedge reverses the %produces the munis 
sign in front of the graph.} 
%%%
(L${}\prec{}$R) of the arrowhead vertices for the two edges issued from the respective wedge top, i.e.\ from the arrowtail.\footnote{\label{FootSunflower}
%\textbf{Example {\textmd{(\cite{Ascona96,skew23})}}.} For $m=1$ and $n=3$, 
The sunflower (\cite{Ascona96,skew23}) is 
$X^{\gamma_3} = (\mathsf{0},\mathsf{1};\mathsf{1},\mathsf{3};\mathsf{1},\mathsf{2}) + 2\cdot (\mathsf{0},\mathsf{2};\mathsf{1},\mathsf{3};\mathsf{1},\mathsf{2})$, here $(m,n)=(1,3)$.%
}
%%%
Each edge carries its own summation index ranging $[1,\ldots,\dim M^d]$ for the 
affine\footnote{\label{FootWhyAffine}
The Poisson manifold is %must be 
affine to make %the 
formulas, differential\/-\/polynomial in Poisson bracket's components and encoded by Kontsevich's graphs, independent of coordinate changes%, now
~$\bx = A\mathbf{x} + \boldsymbol{b}$.}
%%%
Poisson manifold at hand.
An edge decorated by index~$\mathsf{i}$ encodes the derivative $\partial/\partial x^i$ w.r.t.\ the affine coordinates $x^1$,\ $\ldots$,\ $x^d$ of a base point $\bx\in M^d$.
The formula\footnote{\label{Foot NoDimInMKgraph}
Unlike its formula, the Kontsevich graph, with a copy of Poisson bi\/-\/vector $P=P^{ij}\,\partial_{i}\otimes\partial_{j}$ in each internal vertex, does not depend on the dimension~$d$; its topology and the ordering L${}\prec{}$R of outgoing edges in every wedge encodes the operator's formula for every $d\geqslant 2$.}
%%%
of polydifferential operator is the sum (over all indices running $\mathsf{1}$,\ $\ldots$,\ $\mathsf{d}$) of the product of vertices' content.

By the number $d-2$ of its Casimirs $a^k$, the Nambu\/--\/Poisson bracket in Eq.~\eqref{EqNambuBr} is specific to dimension~$d$. Take a microscope and magnify each internal vertex of Kontsevich's graph; it telescopes to a brush: the Levi\/-\/Civita arrowtail vertex contains $\varrho(\bx)\cdot\veps^{i_1\ldots i_d}$, with $\partial_{i_1}\wedge\ldots\wedge\partial_{i_d}$ on the wedge\/-\/ordered $d$-\/tuple of outgoing edges, and (still within the magnified zone) there are $d-2$ terminal vertices with the Casimirs $a^1$,\ $\ldots$,\ $a^{d-2}$ of the respective copy of Nambu\/--\/Poisson bracket.
   %%% TADPOLES.
By definition, a \emph{Nambu micro\/-\/graph} over~$M^d$ is made of whole copies of  Nambu\/--Poisson bracket~\eqref{EqNambuBr} as building blocks.\footnote{\label{FootNambuMKmicrograph}
Every Kontsevich's graph can be expanded to a sum of Nambu micro\/-\/graphs over given $d\geqslant 3$ by working out the Leibniz rules for arrows which entered the internal vertices. Yet not every linear combination of Nambu micro\/-\/graphs over~$M^d$ is Kontsevich's, and not every digraph built over $m$ sinks from suitably many Levi\/-\/Civita and Casimir vertices is Nambu.}
%%%
A Nambu micro\/-\/graph over~$M^d$ consists of $m$ sinks $\mathsf{0}$,\ $\ldots$,\ $\mathsf{m-1}$ and of $n$ copies of bracket~\eqref{EqNambuBr}, which provide the Levi\/-\/Civita vertices (with $\varrho(\bx)\cdot\smash{\veps^{\vec{\imath}}}$) labelled $\mathsf{m}$,\ $\ldots$,\ $\mathsf{m+n-1}$ and the Casimir vertices,\footnote{\label{FootWhoseCasimir}
The label of a Casimir~$a^k$ differs by $k\cdot n$ from the label of 
its `parent' Levi\/-\/Civita vertex.% in the same copy of Nambu structure.
}
%%%
namely: $\mathsf{m+n}$,\ $\ldots$,\ $\mathsf{m+2n-1}$ with~$a^1$, then $\mathsf{m+2n}$,\ $\ldots$,\ $\mathsf{m+3n-1}$ with~$a^2$, etc.
Every Nambu micro\/-\/graph is now encoded by the ordered $n$-\/tuple of ordered lists (of length~$d$) of arrowhead vertices for the $d$-tuples of edges issued from the Levi\/-\/Civita arrowtail vertices.

\begin{ex}\label{Ex3DvanishSunflower}
Over $d=3$, we shall refer from Table~\ref{vanishingtable} to the twelve \emph{vanishing} (as formulas) 1-vector Nambu micro\/-\/graphs build over $m=1$ sink of $n=3$ tridents: the sink is~$\mathsf{0}$, the Levi\/-\/Civita trident tops are $\mathsf{1},\mathsf{2},\mathsf{3}$, and the Casimir vertices $\mathsf{4},\mathsf{5},\mathsf{6}$ (with~$a^1$) are terminal; their encodings are given 
in\cite[Lemma~2]{MSJB}, typeset in boldface.
Each encoding is the ordered list of $n=3$ ordered $(d=3)$-tuples of arrowhead vertices for the triples of edges issued from the arrowtails~$\mathsf{1},\mathsf{2},\mathsf{3}$.
\end{ex}
   %Nos.38,41 are ZERO.

\section{Vanishing micro\/-\/graphs and embedding: $d \hookrightarrow d+1$}\label{SecVanish3D4D}\noindent%
We now see two ways how, from a Nambu micro\/-\/graph over dimension~$d$, we can construct (a linear combination of) Nambu micro\/-\/graph(s) over the next dimension $d+1$. Namely, equip every Levi\/-\/Civita vertex with one extra Casimir $a^{d-1}$ in its new terminal vertex (to which the new, ordered \emph{last} arrow is now sent from its `parent' Levi\/-\/Civita vertex at the top of the $(d+1)$-brush).\\[0.5pt]
%%%
\mbox{ }\quad\textbf{Descendants.}\quad
If, in dimension~$d$, an external arrow --\,issued from another copy of Nambu\/--\/Poisson bi%\/-\/
vector\,-- acted on a Casimir $a^k$, $1\leqslant k \leqslant d-2$, within a subgraph of Nambu structure, then let the external arrow run --\,via Leibniz rule\footnote{\label{FootEdgeOwnLeibniz}
These Leibniz rules work independently one from another for each incoming edge.}%
%%%
\,-- over its old target and the newly attached Casimir vertex with $a^{d-1}$.
Every expansion of Leibniz rule thus yields a linear combination of $(d+1)$-\/dimensional Nambu micro\/-\/graph \emph{descendants} of the originally taken $d$-dimensional Nambu micro\/-\/graph.\\[0.5pt]
   %%% In reserve: {FootFrom2DtoDim}
\mbox{ }\quad\textbf{Embedding.}\quad
After the new Casimirs $a^{d-1}$ are attached, one per each Levi\/-\/Civita vertex, none of the old edges is re\/-\/directed and 
no Leibniz rules are worked out.\footnote{\label{FootExEmbed}
See~\cite[Example~3]{IV} for an illustration: 
$\Gamma_{d=3} \hookrightarrow \smash{ \widehat{ \Gamma }_{d=4} }$; 
similar is Definition~5 and Propositions~4,5 in the paper 
[\href{https://arxiv.org/abs/2409.18875}{I.}] (\href{https://arxiv.org/abs/2409.18875}{\textrm{arXiv:2409.18875}} [q-alg])
or Definition~3 and Example~2 in %the paper
~[\href{https://arxiv.org/abs/2409.15932}{III.}], which is
\href{https://arxiv.org/abs/2409.15932}{arXiv:2409.15932} [q-alg].}
%%%
This yields the \emph{embedding}
$\Gamma_{d} \hookrightarrow \smash{ \widehat{ \Gamma }_{d+1} }$ of the original (micro) graph $\Gamma$ from dimension $d\geqslant 2$ to the Nambu micro\/-\/graph over dimension~$d+1$.

%\begin{rem}\label{RemCrossTermsInEmbed}
Consider the formula of the object (e.g., 1-vector on~$\BBR^{d+1}$) encoded by the Nambu micro\/-\/graph $\smash{ \widehat{ \Gamma }_{d+1} }$ after embedding $\Gamma_d$ as its subgraph. The new edge(s) to new Casimir(s) carry new indice(s).
When \emph{each} new index equals $d+1$, encoding $\partial/\partial x^{d+1} (a^{d-1})$ in the new vertex, the old formula $\phi(\Gamma_d)$ reappears:
\begin{equation}\label{EqWithCrossTerms}
\phi\bigl(\smash{ \widehat{ \Gamma }_{d+1} }\bigr) = \phi(\Gamma_d) \cdot 
\bigl(\partial a^{d-1} / \partial x^{d+1} \smash{\bigr) ^n} + 
\langle \text{cross\/-\/terms} \rangle,
\end{equation}
$n$~being the number of Levi\/-\/Civita vertices.
The cross\/-\/terms are those where %in which 
at least one new $\partial/\partial x^{d+1}$ acts on the content of old vertices from the subgraph~$\Gamma_d$.
%\end{rem}
We discover a curious %strange %counter\/-\/intuitive 
property of Jacobian determinants in brackets~\eqref{EqNambuBr}: the vanishing $\phi(\Gamma_d)=0$ is preserved by the embedding $\Gamma_{d} \hookrightarrow \smash{ \widehat{ \Gamma }_{d+1} }$, 
so $\phi\bigl(\smash{ \widehat{ \Gamma }_{d+1} }\bigr) = 0$; all the cross\/-\/terms cancel out!
%%%
Exploring %In the context of 
the open problem --\,is the tetrahedral graph cocycle %~$\gamma_3$ 
action on the space of %Nambu class~\eqref{EqNambuBr} of 
Poisson brackets \textit{non}trivial\,?\,-- we %let us 
study the vanishing mechanism(s) and the work of (micro-)\/graph embeddings.\footnote{\label{FootVanish4DHam}%26.APR.2024: F.S.+MSJB; 
\textbf{Example.} For $d=4$ and $(m,n)=(0,2)$, the only vanishing (as formula) Hamiltonian is
$H^{(9)}_{d=4} = [\mathsf{1}$,$\mathsf{2}$,$\mathsf{3}$,$\mathsf{5}$;\ $\mathsf{3}$,$\mathsf{4}$,$\mathsf{5}$,$\mathsf{6}]$ from [\href{https://arxiv.org/abs/2409.15932}{III.}, Lemma~16]: here $\mathsf{1}$,$\mathsf{2}$ are Levi\/-\/Civita vertices, $\mathsf{3}$,$\mathsf{4}$ are Casimirs $a^1$, and $\mathsf{5}$,$\mathsf{6}$ are $a^2$; its embedding still vanishes,
$\phi \smash{ \bigl(\widehat{H}^{(9)}_{d=5}\bigr) } = 0$.}

The tetrahedral flow $\smash{ \dot{P} = Q^{\gamma_3}\bigl(P^{\otimes^4}\bigr)}$ needs $1$-vectors $\smash{ \vec{X}^{\gamma_3}\bigl(P^{\otimes^3}\bigr)}$ to be trivialsed (if, indeed, $Q^{\gamma_3}$ is a Poisson coboundary $[\![ P, \smash{ \vec{X}^{\gamma_3} } ]\!]$); to find $\smash{ \vec{X}^{\gamma_3}_{d\geqslant 2} }$ for the Nambu class, we use %Nambu
micro\/-\/graphs on $m=1$ sink and $n=3$ copies of bracket~\eqref{EqNambuBr} over $\BBR^d$. To make the task for $d+1$ smaller, in~\cite{MSJB} and~\cite{IV} we learned to use the descendants of `sunflower' graph from 2D (see footnote~\ref{FootSunflower} on p.~\pageref{FootSunflower}), adjoining the set of $(d+1)$-\/descendants of \emph{vanishing} $1$-vector Nambu micro\/-\/graphs, which were invisible in %the
lower dimension~$d$.
Specifically for $d=3 \hookrightarrow d+1=4$, consider
%\footnote{\label{FootNoVanishVec2D}
%On $\BBR^2$ there are no vanishing $1$-vectors for $n=3$.} 
%%%
%the %exhaustive
%full list of 
sunflower's twelve vanishing descendants: see Example~\ref{Ex3DvanishSunflower} on p.~\pageref{Ex3DvanishSunflower} and Table~\ref{vanishingtable}.\footnote{\label{FootLegendTableEC}
%\underline
\textbf{Legend \textmd{(to Table~\ref{vanishingtable})}.}\quad %\\ %\quad 
Row~1 (R1) contains the indices %numbers 
(in a list of~48) of vanishing 3D micro\/-\/graph descendants $\Gamma\in\text{Van}_{d=3}$(sunflower) from~2D.
%of the `sunflower' graphs from~2D.
%$\bullet$ Row 1 (R1) contains the indices of the micro-graphs $g\in\text{Van}_{d=3}$(sunflower). \\ 
\quad%\qquad %\\
%%%
$\bullet$ Micro\/-\/graphs with %which have 
a nontrivial %non-unit 
automorphism group are labelled by \textbf{aut}.\quad %\\
$\bullet$ Micro\/-\/graphs which %are 
equal %to 
minus themselves under some automorphism
are labelled by \textbf{zero}.\quad %\qquad %\\ 
%%%
$\bullet$ Row~2 (R2) 
contains the number of 4D-descendants of each micro\/-\/graph $\Gamma\in\text{Van}_{d=3}$(sunflower).\quad %\\
%%%
$\bullet$ Row~3 (R3) 
contains the number of vanishing 4D-descendants of each micro\/-\/graph $\Gamma\in\text{Van}_{d=3}$(sunflower).\quad %\qquad %\\ 
%%%
$\bullet$ Row~4 (R4) states %reports %(R4) denotes 
the nature of %the 
vanishing 4D-descendants for each micro\/-\/graph $\Gamma\in\text{Van}_{d=3}$(sunflower);
% $\bullet$ 
the embedding map is denoted by~\textbf{e}%$e$
, the contra\/-\/embedding \textbf{c} %$c$ 
then amounts to $a^1\rightleftarrows a^2$.%
%exclusively redirects arrows to the new Casimirs~$a^2$.
%\\ 
}
%%%%%%%%%%%%%%%%%%%%%%%%%%%%%%%%%%%
\begin{table}[htbp]
\caption[qqqq]
{Vanishing 4D-descendants of 3D micro\/-\/graphs Van$_{d=3}$(sunflower).}\label{vanishingtable}\smallskip
\begin{center}\small
\tabcolsep=3.6pt
\begin{tabular}{l|l|l|l|l|l|l|l|l|l|l|l|l|c}
\hline\hline
%$\dfrac{\text{10}}{\text{aut}}$ & 13 & 20 & 21 & $\dfrac{\text{24}}{\text{aut}}$ & 25 & 29 & $\dfrac{\text{32}}{\text{aut}}$ & $\dfrac{\text{33}}{\text{aut}}$ & 37 & $\dfrac{\text{38}}{\text{zero}}$ & $\dfrac{\text{42}}{\text{zero}}$ 
R1: & 10 & 13 & 20 & 21 & 24 & 25 & 29 & 32 & 33 & 37 & 38 & 41 %42 
& Total \\
 & aut & & & & aut & & & aut & aut & & zero & zero & \mbox{ } \\
\hline
R2: &8   & 2   & 4       & 4       & 8       & 8   & 4       & 8   & 8   & 8       & 16  & 32       & 118               
\\ \hline
%%%
R3: & 2   & 2   & 4       & 4       & 8       & 2   & 4       & 2   & 2   & 8       & 2   & 12       & 54                
\\ \hline
%%%
R4: & e,c & e,c & e,c & e,c & e,c & e,c & e,c & e,c & e,c & e,c & e,c & e,c & e,c \\
 &    &     & ${}+2$ & ${}+2$ & ${}+6$  &     & ${}+2$  &     &     & ${}+6$  &     & ${}+10$ & \textbf{vanish} \\
%%%
%\dfrac{\text{e}}{\text{c}} & \dfrac{\text{e}}{\text{c}} & \dfrac{\text{e}}{\text{c}} +2 & \dfrac{\text{e}}{\text{c}} +2 & \dfrac{\text{e}}{\text{c}} +6 & \dfrac{\text{e}}{\text{c}} & \dfrac{\text{e}}{\text{c}} +2 & \dfrac{\text{e}}{\text{c}} & \dfrac{\text{e}}{\text{c}} & \dfrac{\text{e}}{\text{c}} +6 & \dfrac{\text{e}}{\text{c}} & \dfrac{\text{e}}{\text{c}} +10
%%%
%\thead{e,c} & \thead{e,c} & \thead{e,c} & \thead{e,c \\+ 2} & \thead{e,c \\+ 2} & \thead{e,c \\+ 6} & \thead{e,c} & \thead{e,c \\+ 2} & \thead{e,c} & \thead{e,c} & \thead{e,c \\+ 6} & \thead{e,c} & \thead{e,c \\+ 10} & \thead{\textbf{e,c} \\\textbf{always} \\\textbf{vanish}} 
%%%
\hline
\hline
\end{tabular}\\ %\label{t01}
\end{center}
\normalsize
%\mbox{ }\\[-25pt]\mbox{ }
%\medskip
%%%
\end{table}
%%%%%%%%%%%%%%%%%%%%%%%%%%%%%%%%%%%%%%%%
Of them, Nos.\,38 and 41 are \emph{zero}: they have a symmetry $g\in\Aut(\Gamma)$ that acts on the vertices and edges (thus relabelling the summation indices) and makes $\Gamma=-\Gamma$ w.r.t.\ the wedge ordering L${}\prec{}$M${}\prec{}$R of edges in the tridents, so $\phi(\Gamma)=-\phi(\Gamma)=0$.
The other ten descendants have no symmetry with such effect; being nonzero, they vanish in another way.% We now outline what it~is.%that mechanism~is.
\\[0.5pt]
\mbox{ }\quad\textbf{Case $\Aut(\Gamma)\neq\{\mathbf{1}\}$.}\quad
Four nonzero micro\/-\/graphs in Table~\ref{vanishingtable} still have nontrivial symmetry groups (Nos.\,10,\ 24,\ 32,\ 33). 
We detect\footnote{\label{FootHowVanishByAut}
This claim is verified by brute force calculation, see \cite[Example~5]{IV} for graph No.\,10 in Table~\ref{vanishingtable}.}
%%%
that for each of them, its formula splits into mutually cancelling disjoint pairs of terms; these pairs are marked by exactly those Casimirs which are effectively moved by at least one element $g\neq\mathbf{1}$ from $\Aut(\Gamma)$. In brief, automorphisms %serve as pointers 
highlight %ing 
the key factors.\\[0.5pt]
\mbox{ }\quad\textbf{Case $\Aut(\Gamma) = \{\mathbf{1}\}$.}\quad
For the six remaining nonzero vanishing descendants of the sunflower, the identity $\phi(\Gamma)=0$ is still due to the cancellation of disjoint pairs of terms, without accumulation of longer linear combinations; yet no vertices are marked by the effective action of a symmetry. To reveal this mechanism in full is a standing problem; the impact of topology (in~$\Gamma$) on arithmetic (in $\phi(\Gamma)$) will make possible the study of $d \gtrsim 5$, so far costly.
   %cf.[skew23] vs [II,IV].
   
We argue that these two mechanisms of $\phi(\Gamma_d)=0$ for (non)zero (micro-)\/graphs persist under $\Gamma_d \hookrightarrow \smash{\widehat{\Gamma}_{d+1}}$. For zero $\Gamma\cong %\simeq 
g(\Gamma)=-\Gamma$ the extension of~$g$ after new Casimirs are adjoined is verbatim. For nonzero $\Gamma_d$, %the 
disjoint pairs of terms cancelled out as the summation indices ran, independently one from another, up to $d$. Yet the adjoining of new Casimirs and new summation index in every Levi\/-\/Civita symbol $\veps^{i_1\ldots i_{d+1}}$ does not alter the (grouping of) previously existing factors and terms, only lifting the summation limit to $d+1$. The cancellations work as before, now on a larger set of indices and %on 
a longer range of each index.\footnote{\label{FootExNo10EmbedVanish}
\textbf{Example.}
Taking micro\/-\/graph No.\,10 from Table~\ref{vanishingtable} in $d=3$, one easily upgrades $\phi(\Gamma_{d=3})=0$ to $\phi \smash{\bigl( \widehat{\Gamma}_{d=4} \bigr)} = 0$ by three new factors.}

Let us compare the set $\Van_{d=4}$ of vanishing 
4D-descendants\footnote{\label{FootCountSunflower4D}
There are 324 4D-descendants of the `sunflower'; 54 of them vanish (see~\cite{MSJB}).}
%%%
of the `sunflower' in 2D with, on the other hand, the set of vanishing 4D-descendants of the twelve (from Table~\ref{vanishingtable}) \emph{vanishing} 3D-descendants $\Gamma\in\Van_{d=3}$ of the `sunflower'.\\[2pt]
%%%
\textbf{Proposition~\textbf{1}.}\quad
$\Van_{d=4} = \text{4D-\/descendants}\,\smash{ \bigl( \Van_{d=3} \bigr) }$.\\[2pt]
%%%
This means that in $d+1=4$, we do not seek for the vanishing $1$-vectors anywhere \emph{else}; there appear no new (%detectably 
starting to work) mechanisms of $\phi(\Gamma_{d+1})=0$ w.r.t.\ $\phi(\Gamma_{d})=0$.
This is important (see \cite{IV}): we reduce the intractable problem in dimension $d+1$ before trying to solve it; e.g., $d=5$ is beyond the power of H\'abr\'ok computing cluster at hand.

\textbf{Open problems}, specific to Kontsevich's graph cocycle flows on the (sub)spaces  of Nambu\/--\/Poisson brackets~\eqref{EqNambuBr} on $\BBR^{d\geqslant 3}$, are in particular these:\\[0.5pt]
\mbox{ }$\bullet$\quad If $\smash{\dot{P} = Q^{\gamma_3}_{d+1} \bigl( P^{\otimes^4} \bigr) }$ is trivial, $\smash{Q^{\gamma_3}_{d+1} = [\![ P, \vec{X}^{\gamma_3}_{d+1} \bigl( P^{\otimes^3} \bigr) ]\!] }$, then is $\smash{ \vec{X}^{\gamma_3}_{d+1} }$ found over linear combinations of the `sunflower' descendants from~2D\,?\\[0.25pt]
%%%
\mbox{ }$\bullet$\quad If yes, is $\smash{ \vec{X}^{\gamma_3}_{d+1} }$ found over the union of %micro\/-\/graphs from 
$(d+1)$-\/descendants for the underlying %previous 
solution $%\smash
{ \vec{X}^{\gamma_3}_{d} }$ and 
vanishing $1$-vector micro\/-\/graphs over~$\BBR^d$\,?\\[0.25pt]
%%%
\mbox{ }$\bullet$\quad Do these indispensable (cf.~\cite{IV}) vanishing micro\/-\/graphs stem \emph{only} from the underlying vanishing set in~$\BBR^{d-1}$ through their descendants\,?\\[0.25pt]
%%%
\mbox{ }$\bullet$\quad What is the group\/-\/theoretic and topological mechanism of vanishing for Nambu micro\/-\/graphs and their linear combinations\,?

%\subsection*{Acknowledgements}
\noindent\textbf{Acknowledgements.}\quad
The authors thank the organizers of the XIII International symposium on Quantum Theory and Symmetries -- QTS13 (Yerevan, Armenia, 28~July -- 1~August 2025) for a 
%dynamic and welcoming 
warm atmosphere. 
The authors thank the Center for Information Technology of the University of
Groningen for access to the High Performance Computing cluster, H\'abr\'ok. 
The authors thank F.~Schipper and R.~Buring for helpful discussions and advice.%
%, as well as help with coding in \textsf{gcaops} and optimising the code.  
\footnote{%Participation of 
M.\,S.\,Jagoe Brown was supported at QTS13 by the Master's Research Project funds at the Bernoulli Institute, University of Groningen; %that of 
A.\,V.\,Kiselev was supported by project~135110.}

\appendix
\section{Resilience of the graph calculus in the dimensional shift $d\mapsto d+1$.
Examples}\label{SecExamplesIV}
\noindent%
In a series of examples we examine the mechanism behind vanishing graphs -- graphs whose formulas obtained via the graph calculus are equal to zero.
We detect %propose 
a partial answer behind the vanishing mechanism for a certain class of graphs; and we observe a consistent pattern in how the components of the graph formulas cancel out. 

\begin{ex}[Embedding of 3D graph into 4D]\label{ExEmbed3D4D}
We take the graph built of 3D Nambu--Poisson structures given by the encoding 
$$
e=(0,2,4;1,3,5;1,2,6)
$$
and embed it into 4D (that is, we apply the embedding map to $e$): 
$$
\text{embedding}(e)=(0,2,4,\mathbf{7};1,3,5,\mathbf{8};1,2,6,\mathbf{9}),
$$ where the new Casimirs $a^2\in\{7,8,9\}$ which appear in the 4D Nambu--Poisson structure are in bold font. Recall that each tuple (separated by a semi-colon) in the encoding $e$ corresponds to the outgoing arrows of each Nambu--Poisson structure. %Indeed, we see that 
In embedding($e$) the arrows from the graph built of 3D Nambu--Poisson structures (the first three vertex numbers in each tuple of the encoding $e$) remain as they were, with the only difference being that each structure has an outgoing edge to the new Casimir acquired in the dimensional step 3D$\mapsto$4D (the last vertex number in each tuple). \hfill$\blacksquare$
\end{ex}

%This fact initially seems remarkable due to the mechanism of the graph calculus. That is, we cannot view the embeddings of 3D vanishing sunflower micro-graphs as graphs containing vanishing sub-structures. 
%That is, we cannot guarantee 
It is \textit{a priori} not obvious that %these 
the embeddings of vanishing (micro-)graphs 
%graphs will 
will again vanish due to the vanishing of their sub-structures. 
Indeed, the assembly of formulas using the graph calculus implies the creation of a new family of cross-terms.
%, therefore the formulas have no reason to preserve information from the vanishing sub-structure. 

\begin{ex}\label{ExAppearCrossTermsInEmbed}
    Let us take the formula of the 3D vanishing sunflower micro-graph $g$ with index 10 in Table \ref{vanishingtable} given by the encoding $e_g$, and embed it into 4D: $$e_g=(0,1,4;1,6,5;4,5,6),$$
$$\text{embedding}(e_g)=(0,1,4,\textbf{7};1,6,5,\textbf{8};4,5,6,\textbf{9}),$$ with the index of the new Casimir $a^2$ in bold. We write the inert sum of the formula of $g$ in 3D:$$\phi(g)=\sum_{\vec{\imath},\vec{\jmath},\vec{k}}^{d=3}\varepsilon^{i_1i_2i_3}\varepsilon^{j_1j_2j_3}\varepsilon^{k_1k_2k_3}\varrho^2\varrho_{i_2j_1}a_{i_3k_1}a_{j_3k_2}a_{j_2k_3}\partial_{i_1}(),$$ and the inert sum of the embedding of $g$ into 4D: $$\phi\bigl(\text{embedding}(g)\bigr) =$$ $$\sum_{\vec{\imath},\vec{\jmath},\vec{k}}^{\boldsymbol{d=4}}\varepsilon^{i_1i_2i_3\boldsymbol{i_4}}\varepsilon^{j_1j_2j_3\boldsymbol{j_4}}\varepsilon^{k_1k_2k_3\boldsymbol{k_4}}\varrho^2\varrho_{i_2j_1}a^1_{i_3k_1}a^1_{j_3k_2}a^1_{j_2k_3}\boldsymbol{a^2_{i_4}a^2_{j_4}a^2_{k_4}}\partial_{i_1}(),$$ where the terms concerning the new Casimir $a^2$ and dimension 4D are in bold. That is, each Nambu--Poisson structure which composes the graph embedding$(g)$ in 4D has four outgoing edges (instead of three, as in 3D). Therefore, the indices in the inert sum which correspond to the outgoing edges of the Nambu--Poisson structures will run over $\{1,2,3,4\}$, which will create cross-terms in such a way that we lose track of the formula of the 3D vanishing micro-graph $g$. That is, the 3D formula is reproduced and multiplied by the terms in bold with $\boldsymbol{i_4,j_4,k_4}=4$. But there appear many other terms, when the indices are permuted over $\{1,2,3,4\}$. \hfill$\blacksquare$
\end{ex}

The approach we take %in order 
--\,to investigate why the embeddings of 3D vanishing sunflower micro-graphs vanish\,-- is to understand %investigate 
why they themselves vanish in 3D.
%To this end, we examine the formula of all 3D vanishing sunflower graphs. 

\begin{ex}\label{neatlyvanishing}
    Let us look at the above example of the 3D vanishing sunflower micro-graph $g$ with index 10 in Table \ref{vanishingtable} with a non-trivial automorphism group, where we show in bold the Casimirs on which there acts the non-trivial automorphism group:  
$$\phi(g)=\sum_{\substack{i_1,i_2,i_3\\ j_1,j_2,j_3,\\k_1,k_2,k_3=1}}^{d=3}\varepsilon^{i_1i_2i_3}\varepsilon^{j_1j_2j_3}\varepsilon^{k_1k_2k_3}\varrho^2\varrho_{i_2j_1}a_{i_3k_1}\boldsymbol{a_{j_3k_2}a_{j_2k_3}}\partial_{i_1}().$$
We plug in a certain permutation of the $\vec{\imath}$ terms into the inert sum, $\vec{\imath}=(1,2,3)$, that is $i_1=1,i_2=2,i_3=3$: 

$\bullet \quad \vec{\imath}=(1,2,3)$:$$\sum_{j,k}\varepsilon^{123}\varepsilon^{j_1j_2j_3}\varepsilon^{k_1k_2k_3}\varrho^2\varrho_{2j_1}a_{3k_1}a_{j_3k_2}a_{j_2k_3}\partial_1().$$
We now plug in two consecutive permutations of the $\vec{\jmath}$ terms:

$\bullet\quad(\vec{\imath},\vec{\jmath})=\bigl(\vec{\imath}=(1,2,3),\vec{\jmath}=(1,2,3)\bigr)$:$$\sum_k\varepsilon^{123}\varepsilon^{123}\varepsilon^{k_1k_2k_3}\varrho^2\varrho_{21}a_{3k_1}\boldsymbol{a_{3k_2}a_{2k_3}}\partial_1()=\sum_k\varepsilon^{k_1k_2k_3}\varrho^2\varrho_{21}a_{3k_1}\boldsymbol{a_{3k_2}a_{2k_3}}\partial_1(),$$

$\bullet\quad(\vec{\imath},\vec{\jmath})=\bigl(\vec{\imath}=(1,2,3),\vec{\jmath}=(1,3,2)\bigr)$:$$\sum_k\varepsilon^{123}\varepsilon^{132}\varepsilon^{k_1k_2k_3}\varrho^2\varrho_{21}a_{3k_1}\boldsymbol{a_{2k_2}a_{3k_3}}\partial_1()=\sum_k-\varepsilon^{k_1k_2k_3}\varrho^2\varrho_{21}a_{3k_1}\boldsymbol{a_{2k_2}a_{3k_3}}\partial_1().$$ Here, we can see without expanding the sum any further, that the inert sums of $\bigl(\vec{\imath}=(1,2,3),\vec{\jmath}=(1,2,3)\bigr)$ and $\bigl(\vec{\imath}=(1,2,3),\vec{\jmath}=(1,3,2)\bigr)$ will cancel out due to the Casimirs in bold on which the non-trivial automorphism group acts. Indeed, for a given permutation $\vec{\imath}$ the set of permutations $\vec{\jmath}$
is partitioned into three odd-even pairs, which differ by one
transposition, hence by the +/- sign. These pairs of terms cancel
out for every particular value of the permutation $\vec{k}$.

We find that the other $\bigl(\vec{\imath}=(1,2,3),\vec{\jmath}\bigr)$ terms cancel in the same way for any $\vec{\jmath}$, meaning that all terms with $\vec{\imath}=(1,2,3)$ vanish on their own. There is no cross-cancellation with other $\vec{\imath}$ terms, that is: \begin{align*}\phi(g)&=\sum_{\vec{\imath}=(1,2,3),\vec{\jmath},\vec{k}}^{d=3}+\sum_{\vec{\imath}=(1,3,2),\vec{\jmath},\vec{k}}^{d=3}+\sum_{\vec{\imath}=(2,1,3),\vec{\jmath},\vec{k}}^{d=3}+\sum_{\vec{\imath}=(2,3,1),\vec{\jmath},\vec{k}}^{d=3}+\sum_{\vec{\imath}=(3,2,1),\vec{\jmath},\vec{k}}^{d=3}+\sum_{\vec{\imath}=(3,1,2),\vec{\jmath},\vec{k}}^{d=3}\\&=0+0+0+0+0+0\\&=0.\end{align*} In this example, we showed how the non-trivial automorphism group of the micro-graph acting on the Casimirs induced the formula of the micro-graph to vanish. \hfill$\blacksquare$
\end{ex}

%\noindent \textbf{Why do the embeddings into 4D of $\text{Van}_{d=3}(\text{sunflower})$ vanish?} Now knowing \emph{how} the elements $g\in\text{Van}_{d=3}(\text{sunflower})$ vanish, when we revisit the formula of the embedding of an element $g\in\text{Van}_{d=3}(\text{sunflower})$ into 4D, it \emph{makes sense} why the 4D formula vanishes. 

Let us again look at the same graph as in Example~\ref{neatlyvanishing}: 
$$\phi(g)=\sum_{\vec{\imath},\vec{\jmath},\vec{k}}^{d=3}\varepsilon^{i_1i_2i_3}\varepsilon^{j_1j_2j_3}\varepsilon^{k_1k_2k_3}\varrho^2\varrho_{i_2j_1}a_{i_3k_1}a_{j_3k_2}a_{j_2k_3}\partial_{i_1}(),$$ $$\phi\bigl(\text{embedding}(g)\bigr)=$$ $$ \sum_{\vec{\imath},\vec{\jmath},\vec{k}}^{\boldsymbol{d=4}}\varepsilon^{i_1i_2i_3\boldsymbol{i_4}}\varepsilon^{j_1j_2j_3\boldsymbol{j_4}}\varepsilon^{k_1k_2k_3\boldsymbol{k_4}}\varrho^2\varrho_{i_2j_1}a^1_{i_3k_1}a^1_{j_3k_2}a^1_{j_2k_3}\boldsymbol{a^2_{i_4}a^2_{j_4}a^2_{k_4}}\partial_{i_1}().
$$
%We can now understand why $\phi\bigl(\text{embedding}(g)\bigr)$ vanishes. Because the formulas of all $g\in\text{Van}_{d=3}(\text{sunflower})$ vanish neatly, without cancellation by cross-terms, 
%%%
It is clear that the %this %\emph
{cancellation structure} is %will be 
preserved under the embedding. 

\end{document}